\documentclass[12pt]{amsart}
\usepackage{amscd,amssymb}

\newtheorem{thm}{Theorem}[section]

\newtheorem{lemma}[thm]{Lemma}

\theoremstyle{definition}
\newtheorem{defin}[thm]{Definition}

\theoremstyle{remark}

\numberwithin{equation}{section}

\def\QED{\nobreak\quad\ifmmode\roman{Q.E.D.}\else{\rm Q.E.D.}\fi}

\def\R {{\Bbb R}}

\def\Z {{\Bbb Z}}

\def\N{{\Bbb N}}

\def\T{{\Bbb T}}

\newcommand{\ssk}{\vskip 0.2cm}
\newcommand{\sk}{\vskip 0.5cm}
\begin{document}

\title[Generalized Heisenberg groups]
{Generalized Heisenberg groups and Shtern's question}

\author[Megrelishvili]{Michael Megrelishvili}
\address{Department of Mathematics,
Bar-Ilan University, 52900 Ramat-Gan, Israel}
\email{megereli@math.biu.ac.il}
\urladdr{http://www.math.biu.ac.il/$^\sim$megereli}

\date{October 5, 2004}

\thanks{The author thanks to Israel Science Foundation (grant number 4699).}

\centerline {\small{\rm{Dedicated to the 90th birthday anniversary
 of Professor G. Chogoshvili}}}

 \sk
 \sk


\begin{abstract}
Let $H(X):=(\R \times X) \leftthreetimes X^*$ be the generalized
Heisenberg group induced by a normed space $X$. We prove that $X$
and $X^*$ are relatively minimal subgroups of $H(X)$. We show that
the group $G:=H(L_4[0,1])$ is reflexively representable but weakly
continuous unitary representations of $G$ in Hilbert spaces do not
separate points of $G$. This answers a question of A. Shtern \cite{Sh}.
\end{abstract}


\maketitle

\sk
\section{Groups of Heisenberg type}
\sk

Recall that the classical real 3-dimensional Heisenberg group can
be defined as a linear group of the following matrices:
$$\begin{pmatrix}
1 & a & c\\
0 & 1 & b\\
0 & 0 & 1\\
\end{pmatrix}$$
where $a,b,c \in \R$.  This group is isomorphic to the semidirect
product $(\R \times \R) \leftthreetimes \R $ of $\R \times \R$ and
$\R$.

We need a natural generalization (see, for example, \cite{Re, Mi,
Me1}) which is based on biadditive mappings.

\begin{defin}
Let $E, F$ and $A$ be Hausdorff abelian topological groups and $w:
E \times F \to A$ be a continuous biadditive mapping. Denote by
$$
H(w)= (A \times E) \leftthreetimes F
$$
the semidirect product (say, {\it generalized Heisenberg group}
induced by $w$) of $F$ and the group $A \times E$. The resulting
group, as a topological space, is the product $A \times E \times F
$. The group operation is defined by the following rule: for a
pair
$$u_1=(a_1,x_1,f_1), \hskip 0.4cm u_2=(a_2,x_2,f_2)$$
define
$$u_1 \cdot u_2 = (a_1+a_2+f_1(x_2), x_1+x_2,
f_1 +f_2)$$ where, $f_1(x_2)=w(x_2,f_1)$.
\end{defin}

\sk

Then $H(w)$ becomes a two-step nilpotent Hausdorff topological
group.

We will identify $E$ with $\{0_A \} \times E \times \{0_F \}$ and $F$
with $\{0_A\} \times \{0_E\} \times F$.

Elementary computations for the commutator $[u_1,u_2]$ give
$$[u_1,u_2] = u_1u_2u_1^{-1}u_2^{-1}= (f_1(x_2)-f_2(x_1),0_E,0_F).$$

In the case of a normed space $X$ and the canonical bilinear
function $w: X \times X^* \to \R$ we write $H(X)$ instead of
$H(w)$. Clearly, the case of $H(\R^n)$ (induced by the inner
product $w: \R^n \times \R^n \to \R$) gives the classical
2n+1-dimensional Heisenberg group.

\sk

\section{Relatively minimal subgroups}

\sk

 First we need some definitions concerning the minimality
concept.

\begin{defin}
Let $X$ be a subgroup of a Hausdorff topological group $(G,
\tau)$. We say that $X$ is {\it relatively minimal} in $G$ if
every coarser Hausdorff group topology $\sigma \subseteq \tau$ of
$G$ induces on $X$ the original topology. That is, $\sigma|_X =
\tau|_X$.
\end{defin}

Equivalently, $X$ is relatively minimal in $G$ iff every injective
continuous group homomorphism $G \to P$ into a Hausdorff
topological group induces on $X$ a topological embedding. In
particular, every {\it faithful} (that is, injective) weakly
continuous unitary representation $G \to Is(H)$ for an arbitrary
Hilbert space $H$ induces a topological group embedding (say, {\it
topologically faithful} representation) of a relatively minimal
subgroup $X$.

A {\it minimal group} in the sense of Stephenson \cite{St} and
Do\"\i chinov \cite{Do} is just the group $G$ such that $G$ is
relatively minimal in $G$. Recall some natural examples of minimal
groups: the full unitary group $Is(H)$ (with the weak operator
topology), the symmetric topological group $S_X$, $\Z$ with the
$p$-adic topology, the semidirect product $\R^n \leftthreetimes
\R_+$, every connected semisimple Lie group with finite center
(e.g., $SL_n(\R)$). Note that if $G$ is a locally compact abelian
group with the canonical mapping $w: G \times G^* \to \T$ then the
corresponding generalized Heisenberg group $H(w)=(\T \times G)
\leftthreetimes G^*$ is minimal \cite[theorem 2.11]{Me1}. Hence,
every locally compact abelian group is a group retract of a
minimal locally compact group. For more information see \cite{DPS,
Di, Ga, Me1, me-Gmin}.

The following result seems to be new even for the classical
3-dimensional group $H(\R)$.

\begin{thm} \label{relminimal}
The subgroups $X$ and $X^*$ are relatively minimal in the
generalized Heisenberg group $H(X) = (\R \times X) \leftthreetimes
X^*$ for every normed space $X$.
\end{thm}
\begin{proof}
Let $\tau$ be the given topology of $H(X)$ and suppose that $\sigma
\subseteq \tau$ is a coarser Hausdorff group topology on $H(X)$.
Denote by
$X \times X^*$ the Banach space with respect to the norm
$||(x,f)||:= max\{||x||,||f||\}$.
We prove in fact that the map
$$
q: (H(X), \sigma) \to X \times X^*, \hskip 0.3cm (r,x,f) \mapsto (x,f)
$$
is continuous. This will imply that the natural retractions
$(H(X), \sigma) \to X$ and $(H(X), \sigma) \to X^*$ are
continuous. It guarantees the continuity of the identity maps $(X,
\sigma|_X) \to (X, \tau|_X)$ and $(X^*, \sigma|_{X^*}) \to (X^*,
\tau|_{X^*})$. By the inclusion $\sigma \subseteq \tau$ we get
$\sigma|_X=\tau|_X$ and $\sigma|_{X^*}=\tau|_{X^*}$. This will
mean that $X$ and $X^*$ are relatively minimal subgroups in
$H(X)$.

Assuming the contrary there exists a coarser Hausdorff group
topology $\sigma$ on $H(X)$ such that $q: (H(X), \sigma) \to X
\times X^*$ is not continuous. Since $(H(X), \sigma)$ is a
Hausdorff topological group, one can choose a
$\sigma$-neighborhood $V$ of the neutral element ${\bf
0}:=(0,0_X,0_{X^*})$ such that ${\bf 1}:=(1,0_X,0_{X^*}) \neq
[u,v] $ for every $u,v \in V$, where $[u,v]=uvu^{-1}v^{-1}$.

Since the homomorphism $q$ is not $\sigma$-continuous (at ${\bf
0}$) there exists a positive
$\delta$ such that $q(U)$ is not embedded into the ball
$$
B((0_X,0_{X^*}), \delta):=\{(x,f) \in X \times X^*: \hskip 0.2cm
max\{||x||,||f||\} < \delta \}
$$
 for every $\sigma$-neighborhood $U$ of
$\bf 0$. Then (similar to \cite[Lemma 3.5]{Me1}) it follows that
$q(U)$ is norm-unbounded in $X \times X^*$. Indeed, for every $n
\in \N$ choose a $\sigma$-neighborhood $W$ of $\bf 0$ such that
$$
\underbrace{W \cdot W \cdots W}_n \subseteq U.
$$
As we already know, $q(W)$ is not a subset of $B((0_X,0_{X^*}),
\delta)$. Therefore for every $n \in \N$ there exists a triple
$t_n:=(r_n,x_n,f_n) \in W$ such that the pair $(x_n,f_n)=q(t_n)$
satisfies $||(x_n, f_n)||:= max\{||x_n||, ||f_n||\} \geq \delta$
in $X \times X^*$. Then by the definition of the group operation
in $H(X)$ we get $t_n^n = (s_n,nx_n,nf_n)$ (for some $s_n \in
\R$). Thus, $||q(t_n^n)||=||(nx_n,nf_n)|| \geq n\delta$.
Therefore, $q(U)$ is norm-unbounded. Then there exists a sequence
$S:=\{u_n:=(a_n,y_n,\phi_n)\}_{n \in \N}$ in $U$ such that at
least one of the sets $A:=\{y_n\}_{n \in \N}$ and
$B:=\{\phi_n\}_{n \in \N}$ is unbounded. Suppose first that $A$ is
unbounded. The intersection
$$V_{X^*}:=V \cap X^*$$
is a $\sigma|_{X^*}$-neighborhood of $0_{X^*}$ in $X^*$. Clearly,
$\sigma|_{X^*} \subseteq \tau|_{X^*}$ (= the norm topology of
$X^*$). Therefore, $V_{X^*}$ contains a ball $B(0_{X^*},
\varepsilon_0 )$ of $X^*$ for some $\varepsilon_0 >0$. Since $A$
is norm-unbounded and $B(0_{X^*}, \varepsilon_0) \subseteq
V_{X^*}$, Hahn-Banach theorem implies that the set
$$
<A, V_{X^*}> = \{<y_n,f> = f(y_n) : \hskip 0.3cm n \in \N , \hskip
0.2cm  f \in V_{X^*} \}
$$
is unbounded in $\R$. In fact we have $<A, V_{X^*}> = \R$ because
$<,>: X \times X^* \to \R$ is bilinear and
 $cB(0_{X^*}, \varepsilon_0) \subseteq B(0_{X^*},
\varepsilon_0)$ for every $c \in [-1,1]$.

On the other hand, for every $u_n=(a_n,y_n,\phi_n) \in S \subset
V$ and $f=(0,0_{X},f) \in V_{X^*} \subset V$, the corresponding
commutator $[f,u_n]$ is $(f(y_n),0_X, 0_{X^*})$. Hence,
$[V,V]=\{[a,b] : a,b \in V\}$ contains the subgroup $\R \times
\{0_X\} \times \{0_{X^*}\}$. But then ${\bf 1} \in [V,V]$. This
contradicts our assumption.

The case of unbounded $B=\{\phi_n\}_{n \in \N}$ is similar.
Indeed, observe that we have
$$
<V,V>  \supseteq <V_X, B>  \supseteq <B(0_{X}, \varepsilon_0), B>
= \R
$$
for every $B(0_{X}, \varepsilon_0) \subseteq V_X:=V \cap X$. On
the other hand, $[u_n,x]=(\phi_n(x),0_X, 0_{X^*})$ for every $u_n
\in S$ and every $x:=(0,x,0_{X^*}) \in V_X$. As before, this
implies that ${\bf 1} \in [V,V]$.
\end{proof}

Note that the subgroups $X$ and $X^*$ being relatively minimal in $H(X)$
are not however minimal because every abelian minimal group necessarily is precompact
(see for example \cite{DPS}).

\sk

\section{An application: Shtern's Question}

\sk

Let $Is(X)$ be the group of all linear isometries of a Banach
space $X$. Note that by \cite{Me3} the strong and weak operator
topologies on $Is(X)$ coincide for every reflexive $X$. For some
related recent results about infinite-dimensional representations
of general topological groups we refer to \cite{Pe1}.

As an application of our results we prove the following theorem
which answers a question of A. Shtern \cite{Sh} (for a weaker
version see \cite[Theorem 3.1]{Me2} which states that the additive
topological group $L_4[0,1]$ is not embedded into $Is(H)$ for any
Hilbert space $H$).

\begin{thm} \label{main}
There exists a topological group $G$ such that:
\begin{itemize}
\item
[(a)] Weakly continuous unitary representations of $G$ in Hilbert
spaces do not separate points of $G$.
\item
[(b)] $G$ is a topological subgroup of $Is(X)$ for some {\it
reflexive} Banach space $X$, where $Is(X)$ is endowed with the
strong (equivalently, weak) operator topology.
\end{itemize}
\end{thm}

\ssk

We claim that the desired group $G$ is the Heisenberg group of the
canonical bilinear mapping
$$
w: L_4[0,1] \times L_{\frac{4}{3}}[0,1] \to \R
$$
That is,
$$
G:=H(L_4)= (\R \times L_4) \leftthreetimes L_{\frac{4}{3}}.
$$

First we prove assertion (a) of Theorem \ref{main}.

\begin{lemma}
\label{notseparate} Weakly continuous unitary representations of
the generalized Heisenberg group $G:=H(L_4)$ in Hilbert spaces do not separate
points.
\end{lemma}
\begin{proof}
Assuming the contrary holds true, let $\mathcal{F}=\{h_i: G \to
Is(H_i) : \hskip 0.2cm i \in I\}$ be some point separating set of
weakly continuous unitary representations of $G$. Consider the
corresponding $l_2$-sum $H:=\bigoplus_{i\in I} H_i$ of Hilbert
spaces. Passing to the naturally defined $l_2$-sum of
representations we get a weakly continuous representation $h: G
\to Is(H)$ which is faithful because $\mathcal{F}$ separates the
points.
 By the {\it relative minimality} of the subgroup $L_4$ in $G$
(Theorem \ref{relminimal}), we must conclude that the restriction
map $h|_{L_4}: L_4 \to Is(H)$ necessarily is a topological group
embedding. Therefore, there exists a {\it topologically faithful}
unitary representation of $L_4$ into the unitary group $Is(H)$).
This contradicts \cite[Theorem 3.1]{Me2}.
\end{proof}

\ssk

We say that a map $F: A\times B \rightarrow \R$ has the {\it
Double Limit Property} (in short: DLP) if for every pair of
sequences $\{a_n\}_{n \in \N}, \{b_m\}_{m \in \N}$ in $A$ and $B$
respectively,
$$\lim_m \lim_n F(a_n, b_m)=\lim_n \lim_m F(a_n, b_m)$$
holds whenever both of the limits exist. Let us say that a
continuous function $\phi: G \to \R$ on a topological group $G$
has the DLP if the induced map
$$
F: G \times G \to \R, \hskip 0.2cm F(x,y):=\phi(xy)
$$
has the DLP.

\ssk

 We collect here some auxiliary facts.

\begin{lemma} \label{l:DLP}
\begin{enumerate}
\item
({\it Grothendieck's characterization of weakly almost
periodicity}; see for example, \cite{BJM}) A bounded continuous
function $\phi: G \to \R$ is weakly almost periodic (wap, for
short) iff $\phi$ has the DLP.
\item
(\cite{Sh} and \cite[Fact 2.1]{Me2}) A topological group $G$ can
be embedded into $Is(X)$ (endowed with the strong (equivalently,
weak) operator topology) for some reflexive Banach space $X$ iff
the algebra WAP(G) of all wap functions separates the neutral
element from closed subsets of $G$.
\item
Let $\{a_{n,m}\}_{n,m \in \N}$ be a bounded double sequence of
real numbers. Then there exists a double subsequence
$\{a_{i_n,j_m}\}_{n,m \in \N}$ such that the double limits
$c_1:=\lim_m \lim_n a_{i_n,j_m}$ and $c_2:= \lim_n \lim_m
a_{i_n,j_m}$ both exist.
\item
Let $X$ be a reflexive Banach space with the canonical bilinear
mapping
$$
X \times X^* \to \R, \hskip 0.3cm (x,f) \mapsto f(x).
$$
Then for every pair of bounded subsets $A \subset X, F \subset
X^*$ the restriction map $A \times F \to \R$ has the DLP.
\item
\cite[Lemma 3.4]{Me2} The norm in the Banach space $L_{2k}([0,1])$
has the DLP for every $k \in \N$.
\end{enumerate}
\end{lemma}
\begin{proof} (3): There exists a real number $r$ such that $a_{n,m} \in
[-r,r]$ for every $n,m \in \N$. Consider the sequence $\{v_n\}_{n
\in \N}$, where $v_n:=(a_{n,1},a_{n,2}, \cdots )$. Every $v_n$ can
be treated as an element of the compact metrizable space
$[-r,r]^{\N}$. We can choose a subsequence $\{v_{i_n}\}_{n \in
\N}$ which converges to some $v:=(t_1,t_2,\cdots ) \in
[-r,r]^{\N}$. Then $\lim_n a_{i_n,m}=t_m$ holds for every $m \in
\N$. Since $\{t_m\}_{m \in N}$ is bounded we can find a sequence
$j_1 <j_2< \cdots $ and a real number $c_1 \in [-r,r]$ such that
$\lim_m t_{j_m}=c_1$. It follows that $c_1=\lim_m \lim_n
a_{i_n,j_m}$. Now starting with the double sequence
$\{a_{i_n,j_m}\}_{n,m \in \N}$ we obtain by similar arguments
(switching the roles of $n$ and $m$) its double subsequence such
that the second double limit $c_2$ also exists.

(4): Let $f_n$ and $a_m$ be two sequences respectively in $F$ and
$A$ such that $c_1=\lim_m \lim_n f_n(a_m)$ and $c_2=\lim_n \lim_m
f_n(a_m)$. By the reflexivity, the subsets $F$ and $A$ are
relatively weakly compact in $X^*$ and $X$ respectively.
Eberlein-$\mathrm{\check{S}}$mulian theorem implies that $F$ and
$A$ are relatively sequentially compact with respect to the weak
topologies. Hence there exist a subsequence $f_{i_n}$ of $f_n$ and
a subsequence $a_{j_m}$ of $a_m$ such that the weak limits $\lim_n
f_{i_n}=f$ and $\lim_m a_{j_m}=a$ are defined in $X^*$ and $X$.
Since the map $(X,weak) \times (X^*,weak^*) \to \R$ is separately
continuous we obtain
$$
c_1=\lim_m \lim_n f_n(a_m)=\lim_m \lim_n f_{i_n}(a_{j_m})=\lim_m
f(a_{j_m})=f(a)
$$
$$
c_2=\lim_n \lim_m f_n(a_m)=\lim_n \lim_m f_{i_n}(a_{j_m})=\lim_n
f_{i_n}(a)=f(a).
$$
Thus,
$c_1=c_2$, as desired.
\end{proof}

\sk

The norm of the Banach space $c_0$ does not satisfy the DLP.
Indeed, define $u_n:=e_n $ (the standard basis vectors) and
$v_m:=\sum_{i=1}^m e_i$. Then
$$
1=\lim_m \lim_n ||u_n + v_m|| \neq \lim_n \lim_m ||u_n + v_m||=2.
$$

Note also that by \cite{Me:H} there exists a non-trivial Hausdorff
topological group $G$ (namely, the group $G:=H_+[0,1]$ of all
orientation preserving selfhomeomorphisms of the closed interval
$[0,1]$) such that every (weakly) continuous representation $h: G
\to Is(X)$ is constant for every reflexive $X$. Equivalently,
every wap function on $G$ is constant.

\sk

 Now we prove assertion (b) of Theorem
\ref{main}.

\begin{lemma} \label{refl}
The generalized Heisenberg group $G:=H(L_4)$ is
reflexively representable. That is, there exists a reflexive
Banach space $X$ such that $G$ is a topological subgroup of
$Is(X)$ endowed with the strong (equivalently, weak) operator
topology.
\end{lemma}
\begin{proof}
Define the following continuous bounded real valued function
$$
\phi: G=(\R \times L_4) \leftthreetimes L_{\frac{4}{3}} \to \R,
\hskip 0.4cm \phi(r,x,f)=\frac{1}{1+|r|+||x||+||f|| }
$$
This function separates the neutral element ${\bf 0}=(0,0,0)$ from
closed subsets ${\bf 0} \notin A$ in $G$. So, by Lemma
\ref{l:DLP}.2, it suffices to establish that $\phi$ is wap.
Assuming the contrary holds true we get by Lemma \ref{l:DLP}.1
that $\phi$ does not satisfy the DLP. Therefore there exist two
sequences $u_n=(a_n,x_n,p_n), \hskip 0.2cm v_m=(b_m,y_m,q_m)$ in
$G$ such that for the double sequence
$$\phi(u_n v_m)=\frac{1}{1+|a_n+b_m+p_n(y_m)|+||x_n+y_m||+||p_n+q_m||}.$$
the double limits
$$
s_1:=\lim_m \lim_n \phi(u_n v_m),  \hskip 0.3cm s_2:=\lim_n \lim_m
\phi(u_n v_m)
$$
exist but $s_1 \neq s_2$.

We can assume without restricting of generality (up to the
subsequences) that all sequences
$$a_n, x_n, p_n, b_m, y_m, q_m$$ are bounded. Indeed,
if one of the sequences above is unbounded, passing to appropriate
subsequences, we get that the corresponding double limits both are
$0$. Actually we can and do assume that the sequences $a_n$ and
$b_m$ converge in $\R$. Moreover, by Lemma \ref{l:DLP}.3 we can
suppose in addition that each of the following three (bounded)
real valued double sequences
$$
||p_n+q_m||, \hskip 0.3cm  ||x_n+y_m||, \hskip 0.3cm p_n(y_m)
$$
have double limits. Now it suffices to show that in each one of
the cases (a), (b), and (c) below the corresponding double limits
are the same. This will imply that $s_1=s_2$, providing the
desired contradiction.

\begin{itemize}
 \item [(a)] First we check that $ \lim_m \lim_n ||p_n+q_m||=\lim_n
\lim_m ||p_n+q_m||$. We have to show that the function
$L_{\frac{4}{3}} \to \R, f \mapsto ||f||$ has the DLP. It suffices
to establish that the function $f \mapsto e^{-||f||^p}$ has the
DLP on $L_{\frac{4}{3}}$. But this function is positive definite
on $L_p$ for every $1\leq p \le 2$ by a classical result of
Shoenberg \cite{Sho}. On the other hand every continuous positive
definite function on a topological group is wap (see
\cite[Corollary 3.3]{Bu}).
\item [(b)] Observe also that the map $L_4 \to \R,  x \mapsto ||x||$ has
the DLP by Lemma \ref{l:DLP}.5. Therefore, $ \lim_m \lim_n
||x_n+y_m||=\lim_n \lim_m ||x_n+y_m||$ holds.
\item [(c)] Finally we show that $ \lim_m \lim_n p_n(y_m) = \lim_n \lim_m
p_n(y_m).$ Note that, as before, by our assumptions, the sequences
$p_n$ and $y_m$ are bounded in $L_{\frac{4}{3}}$ and $L_4$
respectively and the limits $ c_1:=\lim_m \lim_n p_n(y_m), \hskip
0.2cm c_2:=\lim_n \lim_m p_n(y_m)$ both exist.
 Now the equality $c_1=c_2$ directly follows by Lemma
 \ref{l:DLP}.4.
\end{itemize}
\end{proof}

\ssk

The proof of Theorem \ref{main} follows now from Lemmas
\ref{notseparate} and \ref{refl}.

\ssk

%Acknowledgement: I would like to thank to my friends in A.
%Razmadze Mathematical Institute (Tbilisi, Georgia) where I spent
%many memorable years (1982-1991).

 Finally we thank D. Dikranjan, N. Krupnik and A. Shtern for useful suggestions.

\sk
\sk

\bibliographystyle{amsplain}

\end{document}